\documentclass{elsart3-1}

\usepackage{amssymb}

\usepackage[english,francais]{babel}

\usepackage{hyperref}

\newcommand{\DeclareMathOperator}[2]{\def#1{\mathop{\rm#2}}}
\newcommand{\integer}{{\mathord{\mathbb{Z}}}}

\newcommand{\rational}{{\mathord{\mathbb{Q}}}}
\newcommand{\real}{{\mathord{\mathbb{R}}}}
\newcommand{\complex}{{\mathord{\mathbb{C}}}}

\DeclareMathOperator{\Homeo}{Homeo}

\newcommand{\Qrank}{\mathop{\mbox{\upshape$\rational$-rank}}}
\newcommand{\Rrank}{\mathop{\mbox{\upshape$\real$-rank}}}
\newcommand{\Frank}{\mathop{\mbox{\upshape$\F$-rank}}}

\DeclareMathOperator{\SL}{SL}
\DeclareMathOperator{\PSL}{PSL}
\DeclareMathOperator{\Sp}{Sp}
\newcommand{\ints}{\mathcal{O}}

\newcommand{\F}{F}

\newcommand{\pref}[1]{{\upshape(}\ref{#1}{\upshape)}}
\newcommand{\fullref}[2]{{\ref{#1}\pref{#1-#2}}}

\newcommand{\thmrefer}[1]{\renewcommand\thethm
  {\protect\ref{#1}$'$}\addtocounter{thm}{-1}}

\theoremstyle{definition}
\newtheorem{notation}[thm]{Notation}

\def\bmatrix#1{\left[\matrix{#1}\right]}

\begin{document}

\begin{frontmatter}

\selectlanguage{english}

\title
 {Isotropic nonarchimedean $S$-arithmetic groups are not left orderable}

\vspace{-2.6cm}

\selectlanguage{francais}

\title
 {Groupes $S$-arithm\'etiques non-archim\'ediens isotropes ne sont pas
ordonn\'es \`a gauche}

\selectlanguage{english}

\author{Lucy Lifschitz}

\address{Department of Mathematics,
 University of Oklahoma,
 Norman, Oklahoma, 73019, USA}
 \ead{LLifschitz@math.ou.edu}
 \ead[url]{http\char'72//www.math.ou.edu/\char'176llifschitz/}

\author{Dave Morris}

\address{Department of Mathematics and Computer Science,
 University of Lethbridge,
 Lethbridge, Alberta, T1K~3M4, Canada} \address{Department of Mathematics,
 Oklahoma State University,
 Stillwater, Oklahoma, 74078, USA} \ead{Dave.Morris@uleth.ca}
 \ead[url]{http\char'72//people.uleth.ca/\char'176dave.morris/}

\begin{abstract}
 If $\ints$ is either $\integer[ \sqrt{r} ]$ or $\integer[1/r]$, where $r
> 1$ is any square-free natural number, we show that no finite-index
subgroup of $\SL(2,\ints)$ is left orderable. (Equivalently, these
subgroups have no nontrivial orientation-preserving actions on the real
line.) This implies that if $G$ is an isotropic $\F$-simple algebraic
group over an algebraic number field~$\F$, then no nonarchimedean
$S$-arithmetic subgroup of~$G$ is left orderable. Our proofs are based on
the fact, proved by B.~Liehl, that every element of $\SL(2,\ints)$ is a
product of a bounded number of elementary matrices.

\vskip 0.5\baselineskip

\selectlanguage{francais}


\noindent{\bf R\'esum\'e}

\vskip 0.5\baselineskip

\noindent
 Si $\ints$ est soit $\integer[ \sqrt{r} ]$ ou soit $\integer[1/r]$, o\`u
$r > 1$ est un entier positif sans carr\'e, nous prouvons qu'aucun
sous-groupe d'indice fini de $\SL(2,\ints)$ n'est ordonn\'e \`a gauche.
(En d'autres mots, les sous-groupes d'indice fini de $\SL(2,\ints)$ ne
poss\`edent pas d'action non triviale sur la droite respectant
l'orientation.) Cela implique que si $G$ est un groupe alg\'ebrique
$F$-simple isotrope, d\'efini sur un corps de nombres~$F$, alors aucun
sous-groupe $S$-arithm\'etique non-archim\'edien de $G$ n'est ordonn\'e
\`a gauche. La d\'emonstration est fond\'ee sur le fait, due \`a
B.~Liehl, que chaque \'el\'ement de $\SL(2,\ints)$ est le produit d'un
nombre born\'e de matrices \'el\'ementaires.
 \end{abstract}

\end{frontmatter}

\selectlanguage{english}

\section{Introduction}

It is known \cite{Witte-AOC} that finite-index subgroups of
$\SL(3,\integer)$ or $\Sp(4,\integer)$ are not left orderable. (That is,
there does not exist a total order~$\prec$ on any finite-index
subgroup, such that $ab \prec ac$ whenever $b \prec c$.) More generally,
if $G$ is a $\rational$-simple algebraic $\rational$-group, with $\Qrank
G \ge 2$, then no finite-index subgroup of $G_{\integer}$ is left
orderable. It has been conjectured that the restriction on $\Qrank$ can
be replaced with the same restriction on $\Rrank$, which is a much weaker
hypothesis:

\begin{conj} \label{GZnotLOConj}
 If $G$ is a $\rational$-simple algebraic $\rational$-group, with $\Rrank
G \ge 2$, then no finite-index subgroup~$\Gamma$ of $G_{\integer}$ is left
orderable.
 \end{conj}

It is natural to propose an analogous conjecture that replaces $\integer$
with a ring of $S$-integers, and has no restriction on the $\Rrank$:

\begin{conj} \label{GZSnotLOConj}
 If $G$ is a $\rational$-simple algebraic $\rational$-group, and
$\{p_1,\ldots,p_n\}$ is any nonempty set of prime numbers, then no
finite-index subgroup~$\Gamma$ of $G_{\integer[1/p_1,\ldots,1/p_n]}$ is
left orderable.
 \end{conj}

We prove Conjecture~\ref{GZSnotLOConj} under the additional assumption
that $\Qrank G \ge 1$:

\begin{thm} \label{GZSnotLOThm}
 If $G$ is a $\rational$-simple algebraic $\rational$-group, with $\Qrank
G \ge 1$, and $\{p_1,\ldots,p_n\}$ is any nonempty set of prime numbers,
then no finite-index subgroup~$\Gamma$ of
$G_{\integer[1/p_1,\ldots,1/p_n]}$ is left orderable.

More generally, if $G$ is an $\F$-simple algebraic group over
an algebraic number field~$\F$, with $\Frank G \ge 1$, then no
nonarchimedean $S$-arithmetic subgroup~$\Gamma$ of~$G$ is left orderable.
 \end{thm}

We also prove some cases of Conjecture~\ref{GZnotLOConj} (with $\Qrank G =
1$). For example, we consider $\rational$-forms of $\SL(2,\real) \times
\SL(2,\real)$:

\begin{thm} \label{SL2archNotLO}
 If $r > 1$ is any square-free natural number, then no finite-index
subgroup~$\Gamma$ of $\SL \bigl( 2,\integer[ \sqrt{r} ] \bigr)$ is left
orderable.
 \end{thm}

In geometric terms, the theorems can be restated as the nonexistence of
orientation-preserving actions on the line:

 \begin{cor}
 If $\Gamma$ is as described in Theorem~\ref{GZSnotLOThm} or
Theorem~\ref{SL2archNotLO}, then there does not exist any nontrivial
homomorphism $\varphi \colon \Gamma \to \Homeo^+(\real)$.
 \end{cor}

Combining this corollary with an important theorem of \'E.~Ghys
\cite{Ghys-reseaux} yields the conclusion that every
orientation-preserving action of~$\Gamma$ on the circle~$S^1$ is of an
obvious type; any such action is either 
 virtually trivial
 or
 semiconjugate to an action by linear-fractional transformations,
obtained from a composition $\Gamma \rightarrow \PSL(2,\real)
\hookrightarrow \Homeo^+(S^1)$.
  See \cite{Ghys-AOCsurvey} for a discussion of the general
topic of group actions on the circle.

It has recently been proved that certain individual arithmetic groups are
not left orderable (see, e.g., \cite{notRO}), but our results
apparently provide the first new examples in more than ten years of
arithmetic groups that have no left-orderable subgroups of finite index.
They are also the only known such examples that have $\Qrank 1$.

The theorems are obtained by reducing to the fact, proved by B.~Liehl
\cite{Liehl}, that if $\ints = \integer[1/(p_1\ldots p_n)]$ or
$\ints = \integer[\sqrt{r}]$, then $\SL(2,\ints)$ has bounded generation
by unipotent elements. (That is, the fact that $\SL(2,\ints)$ is the
product of finitely many of its unipotent subgroups. For the
general case of Theorem~\ref{GZSnotLOThm}, we also note that $\Gamma$
contains a finite-index subgroup of $\SL \bigl( 2,\integer[1/p] \bigr)$,
for some prime~$p$.)  We are able to prove the same reduction for certain
other groups:

\begin{thm} \label{BddUOrbit}
 Suppose $\Gamma$ is a finite-index subgroup of either
 \begin{enumerate}
 \item \label{BddUOrbit-nonarch}
 $\SL \bigl( 2, \integer[1/r] \bigr)$, for some natural number $r > 1$,
 or
 \item \label{BddUOrbit-arch}
 $\SL( 2,
\ints )$, where $\ints$ is the ring of integers of a number
field~$\F$, and $\F$ is neither $\rational$ nor an imaginary quadratic
extension of~$\rational$,
 or
 \item \label{BddUOrbit-SL3}
 an arithmetic subgroup of a quasi-split $\rational$-form of the
$\real$-algebraic group $\SL(3,\real)$.
 \end{enumerate}
 If $\varphi \colon \Gamma \to \Homeo^+(\real)$ is any
homomorphism, and $U$ is any unipotent subgroup of~$\Gamma$, then every
$\varphi(U)$-orbit on~$\real$ is bounded.
 \end{thm}

\begin{cor} Suppose 
 \begin{itemize}
 \item $\Gamma$ is as described in Thm.~\ref{BddUOrbit}, and 
 \item $\Gamma$ is commensurable to a group that has bounded generation by
unipotent elements. 
 \end{itemize}
 Then every homomorphism $\varphi \colon \Gamma \to \Homeo^+(\real)$ is
trivial. Therefore, $\Gamma$ is not left orderable.
 \end{cor}

Assuming a certain generalized Riemann Hypothesis, G.~Cooke and
P.~J.~Weinberger \cite{CookeWeinberger} proved that the groups described
in part \pref{BddUOrbit-arch} of Thm.~\ref{BddUOrbit} do have bounded
generation by unipotent elements. Thus, if this generalized Riemann
Hypothesis holds, then finite-index subgroups of these groups are not
left orderable. See \cite{Liehl} for relevant results on bounded
generation that do not rely on any unproved hypotheses, and see
\cite{Murty} for a recent discussion of bounded generation.

\section{Proof of Theorem~\fullref{BddUOrbit}{nonarch}}

\begin{notation}
 For convenience, let
 $$ \mbox{$\overline{u} = \bmatrix{1 & u \cr 0 & 1}$,
 \qquad $\underline{v} = \bmatrix{1 & 0 \cr v & 1}$,
 \qquad $\hat{s} = \bmatrix{s & 0 \cr 0 & 1/s}$
 }$$
  for $u,v \in \integer[1/r]$ and $s \in \{\, r^n \mid n \in
\integer\,\}$.
 \end{notation}

Suppose some $\varphi(U)$-orbit on~$\real$ is not bounded above. (This
will lead to a contradiction.) Let us assume $U$ is a maximal unipotent
subgroup of~$\Gamma$.

Let $V$ be a subgroup of~$\Gamma$ that is conjugate to~$U$, but is not
commensurable to~$U$. Then $V_{\rational} \neq U_{\rational}$. Because
$\Qrank \SL(2,\rational) = 1$, this implies that $V_{\rational}$ is
opposite to~$U_{\rational}$. Therefore, after replacing $U$ and~$V$ by a
conjugate under $\SL(2,\rational)$, we may assume
 $$ \mbox{$U = \{\, \overline{u} \mid u \in \integer[1/r] \,\} \cap
\Gamma$
 \qquad and \qquad
 $V = \{\, \underline{v} \mid v \in \integer[1/r] \,\} \cap \Gamma$.}$$

Because $V$ is conjugate to~$U$, we know that some $\varphi(V)$-orbit is
not bounded above. Let 
 $$ x_U = \sup \{\, x \in \real \mid \mbox{the $\varphi(U)$-orbit of~$x$ is
bounded above} \,\} < \infty $$
 and
 $$ x_V = \sup \{\, x \in \real \mid \mbox{the $\varphi(V)$-orbit of~$x$ is
bounded above} \,\} < \infty. $$
 Assume, without loss of generality, that $x_U \ge x_V$.

Fix some $s = r^n > 1$, such that $\hat s \in \Gamma$, and let $B =
\langle \hat s \rangle U$. Because $\langle \hat s \rangle$
normalizes~$U$, this is a subgroup of~$\Gamma$. 
 Note that $\varphi(B)$ fixes $x_U$, so it acts on the interval
$(x_U,\infty)$. Since $\varphi(B)$ is nonabelian, it is well known (see,
e.g., \cite[Thm.~6.10]{Ghys-AOCsurvey}) that some nontrivial element of
$\varphi(B)$ must fix some point of $(x_U,\infty)$. In
 fact, it is not difficult to see that each element of $\varphi(B)
\smallsetminus \varphi(U)$ fixes some point of $(x_U,\infty)$. In
particular, $\varphi(\hat s)$ fixes some point~$x$ of $(x_U,\infty)$. 

The left-ordering of any additive subgroup of~$\rational$ is unique (up
to a sign), so we may assume that 
 $$ \mbox{$\varphi (\overline{u_1}) x <  \varphi (\overline{u_2}) x
\Leftrightarrow u_1 < u_2 $
 \qquad and \qquad
 $ \varphi (\underline{v_1}) x <  \varphi (\underline{v_2}) x
\Leftrightarrow v_1 < v_2 $.} $$
 The $\varphi(U)$-orbit of~$x$ is not bounded above (because $x > x_U$),
so we may fix some $u_0,v_0 > 0$, such that 
 $$ \varphi(\underline{v_0})x < \varphi(\overline{u_0}) x .$$
 For any $\underline{v} \in V$, there is some $k \in \integer^+$, such
that $v < s^{2k} v_0$.
 Then, because $\varphi(\hat s)$ fixes~$x$ and $s^{-2k} < 1$, we have
 \begin{eqnarray*}
 \varphi(\underline{v}) x
 &<& \varphi(\underline{s^{2k}v_0})x
 = \varphi(\hat s^{-k} \underline{v_0} \hat s^{k}) x
 = \varphi(\hat s^{-k}) \varphi( \underline{v_0} )x 
 \cr &<& \varphi(\hat s^{-k}) \varphi(  \overline{u_0})x 
 = \varphi(\hat s^{-k} \overline{u_0}\hat s^{k})x
 = \varphi(\overline{s^{-2k} u_0}) x
 < \varphi(\overline{u_0})x 
 . \end{eqnarray*}

So the $\varphi(V)$-orbit of~$x$ is bounded above by
$\varphi(\overline{u_0})x$. This contradicts the fact that $x > x_U \ge x_V$.

\section{Other parts of Theorem~\ref{BddUOrbit}}

\pref{BddUOrbit-arch} The above proof of Case~\pref{BddUOrbit-nonarch}
needs only minor modifications to be applied with a ring $\ints$ of
algebraic integers in the place of $\integer[1/r]$. (We choose $s =
\omega^n$, where $\omega$ is a unit of infinite order in~$\ints$.) 
 The one substantial difference between the two cases is that the
left-ordering of the additive group of~$\ints$ is far from unique ---
there are infinitely many different orderings. Fortunately, we are
interested only in left-orderings of  $U = \{\,\overline{u} \mid u \in
\ints\,\} \cap \Gamma$ that arise from an unbounded $\varphi(U)$-orbit, and
it turns out that any such left-ordering must be invariant under
conjugation by~$\hat s$. The left-ordering must, therefore, arise from a
field embedding~$\sigma$ of $\F$ in~$\complex$ (such that $\sigma(s)$ is
real whenever $\hat s \in \Gamma$), and there are only finitely many such
embeddings. Hence, we may replace $U$ and~$V$ with two conjugates of~$U$
whose left-orderings come from the same field embedding (and the same
choice of sign).

\medbreak

\pref{BddUOrbit-SL3} A serious difficulty prevents us from applying the
above proof to quasi-split $\rational$-forms of $\SL(3,\real)$. Namely,
the reason we were able to obtain a contradiction is that if
$\overline{u_0}$ is upper triangular, $\underline{v}$ is lower triangular,
$\hat s$ is diagonal, and $\lim_{k \to \infty} \hat s^{-k} \overline{u_0}
\hat s^k = \infty$ under an ordering of~$\Gamma$, then $\lim_{k \to
\infty} \hat s^{-k} \underline{v} \hat s^k = e$. Unfortunately, the
``opposition involution" of $\SL(3,\real)$ causes the calculation to
result in a different conclusion in case~\pref{BddUOrbit-SL3}: if $\hat
s^{-k} \overline{u_0} \hat s^k$ tends to~$\infty$, then $\hat s^{-k}
\underline{v} \hat s^k$ also tends to~$\infty$. Thus, the above simple
argument does not immediately yield a contradiction.

Instead, we employ a lemma of M.~S.~Raghunathan
\cite[Lem.~1.7]{Raghunathan-CSP2} that provides certain nontrivial
relations in~$\Gamma$. These relations involve elements of both~$U$
and~$V$; they provide the crucial tension that leads to a contradiction.

\section*{Acknowledgements}
 The authors would like to thank A.~S.~Rapinchuk for helpful suggestions.
 D.~M.\ was partially supported by a grant from the National Science
Foundation (DMS--0100438).

\end{document}